\documentclass[]{birkmult}

 \theoremstyle{definition}
 
 \theoremstyle{remark}

 \numberwithin{equation}{section}
\begin{document}
\def\di{\displaystyle}
\def\bv{{\mathbf v}}
\title{On the properties of invariants of forms}
\author{Mehdi Nadjafikhah}
\address{School of
Mathematics, Iran University of Science and Technology, Narmak-16,
Tehran, I.R.Iran.}
\email{m\_nadjafikhah@iust.ac.ir}
\author{Parastoo Kabi-Nejad}
\email{parastoo\_ kabinejad@yahoo.com}
\date{ }
\begin{abstract}
This paper is devoted to a discussion of specific properties of
invariants in the theory of forms.

\medskip \noindent {\bf A.M.S. 2000 Subject Classification:} 34C14, 22F30, 19GXX.

\end{abstract}
\maketitle
%
\section{Introduction}
At first, we point out briefly a geometrical interpretation of the
theory of forms in two variable. Clearly, the theory of binary
forms is identical with the geometry of a line (or a bundle of
lines and planes). Analogously, one realizes that a ternary form
is identical to the geometry of the plane, namely of algebraic
curves. for instance, the elliptic curve $y^2=x^3+ax+b$ which
plays a key role in analytic function theory \cite{[1]}, can be
identified as the algebraic variety corresponding to a certian
homogeneous ternary cubic. See \cite{[2], [3], [12]} for the
remarkable connection between elliptic curves, modular forms, and
Fermat's Last theorem.
 Finaly, the theory of quaternary forms
is an essential aim in studying the geometry of space, especially
of algebraic surfaces. The forms with more than three variables do
not admit such a geometric interpretation. \cite{[6]} Theory of
invariants plays a vital role in the theory of forms.It's origins
can be traced back to Cayley (1845),
 who used the term "hyperdeterminants" for functions possessing the invariant property.
 So, we are
led to unearth essentially new and deep properties of forms
through the application of this theory to the theory of forms.
\section{Properties of invariants}
In the following we will initially limit ourselves to the theory
of binary forms and use them to clarify the general concepts. The
generalization of binary forms to forms with arbitrarily many
variables poses no difficulties in most cases.
We always write the general binary form of order n as:
$$f^{(n)}=a_0x^n+a_1x^{n-1}y+\cdots+a_ny^n.$$
\paragraph{Example 2.1}Consider the homogeneous quadratic
polynomial
\begin{eqnarray}
Q(x,y)=ax^2+2bxy+cy^2
\end{eqnarray}
in two variables x,y. Clearly we can recover the inhomogeneous
quadratic polynomial from the assotiated quadratic form and vice
versa.

A glance observation shows that any invertible linear
transformation
\begin{eqnarray}
\overline{x}=\alpha x+\beta y,\qquad \overline{y}=\gamma x+\delta
y,\qquad \alpha \delta -\beta \gamma\neq 0
\end{eqnarray}
will map a homogeneous polynomial in $\overline{x}$ and
$\overline{y}$ according to
\begin{eqnarray}
\overline{Q}(\overline{x},\overline{y})=\overline{Q}(\alpha
x+\beta y,\gamma x+\delta y)=Q(x,y).
\end{eqnarray}
Remarkably, the discriminant of the transformed polynomial is
directly related to that of the original quadratic form. A
straightforward computation verifies that they agree up to the
square of the determinant of the coefficient matrix for the linear
transformation (2):
\begin{eqnarray}
\triangle=ac-b^2={(\alpha \delta -\beta
\gamma)}^2(\overline{a}\overline{c}-{\overline{b}}^2)={(\alpha
\delta -\beta \gamma)}^2\overline{\triangle}
\end{eqnarray}
Therefore, by induction the defining equation of the invariants is
\begin{eqnarray}
I(a_0',a_1',\cdots,a_n')=\delta^pI(a_0,a_1,\cdots,a_n)\nonumber
\end{eqnarray}
\paragraph{Definition 2.2}
Expressions $g=\nu_0+\nu_1+\cdots +\nu_n$ and
$\nu_1+2\nu_2+3\nu_3+\cdots+n\nu_n$ are called the degree and the
weight of a function
$$I(a_0,a_1,\ldots,a_n)=\sum Z_{\nu_0\nu_1\cdots\nu_n}a_0^
{\nu_0}a_1^{\nu_1}\ \cdots a_n^{\nu _n},$$ respectively, where
$Z_{\nu_0 \nu_1 \cdots\nu_n}$ are numerical coefficients. See
\cite{[6]} for more details.
\paragraph{Note:}
In this paper, we shall assume that all terms of
a function have the same weight.
\paragraph{Lemma 2.3}
{\it Every linear transformation of binary forms can be composed
of the following three types of linear transformations :}
\begin{eqnarray}
x&=&\alpha x', \qquad y= \delta y', \\ x&=&x'+\beta y, \qquad y=y',\\
x&=&x',\qquad y=\gamma x'+y'.
\end{eqnarray}
\paragraph{Theorem 4.2}
{\it Every invariant of a form must be homogeneous in the
coefficients and every term must have degree $g=\frac{2p}{n}$,
where $p$ is the exponent of the transformation determinant by
which the invariant changes under substitution of the transformed
coefficients. Furthermore, all terms must have the same weight,
which is also equal to $p$.}

\medskip\noindent{\it Proof:} Since a general linear transformation is
composed of these three types of transformation, application of a
general transformation will in itself not lead to any new results.

Apply the linear transformation
\begin{eqnarray}
x=\alpha x',\qquad\qquad y=\delta y'\nonumber
\end{eqnarray}
to the binary form
$$f(x,y)=a_0x^n+{n\choose 1}a_1x^{n-1}y+\cdots+{n\choose i}a_ix^{n-i}y^i+\cdots+a_ny^n$$
The given form then becomes:
\begin{eqnarray}
f&=&f(\alpha x',\delta y')=a_0\alpha^n
{x'}^n+{n\choose1}a_1\alpha^{n-1}{x'}^{n-1}\delta
y'+\cdots \nonumber \\
&&\cdots+{n\choose
i}a_i\alpha^{n-i}{x'}^{n-i}\delta^i{y'}^i+\cdots a_n\delta^n{y'}^n\nonumber\\
&=&a_0'{x'}^n+{n\choose 1}a_1'x'^{n-1}y'+\cdots +{n\choose
i}a_i'{x'}^{n-i}{y'}^i+\cdots +a_n'{y'}^n,\nonumber
\end{eqnarray}
where
\begin{eqnarray}
a_0'&=&a_0\alpha^n,\nonumber\\
a_1'&=&a_1\alpha^{n-1}\delta,\nonumber\\
&\vdots&\nonumber\\
a_i'&=&a_i\alpha^{n-i}\delta^i,\nonumber\\
&\vdots&\nonumber\\
a_n'&=&a_n\delta^n.\nonumber
\end{eqnarray}
Now, let
$$I(a_0,a_1,\cdots,a_n)=\sum Z_{\nu_0\nu_1\cdots\nu_n}{a_0}^{\nu_0}{a_1}^{\nu_1}\cdots {a_n}^{\nu_n}$$
be an invariant of the form f, where the Z are numerical
coefficients.

Then we have the identity
\begin{eqnarray}
I(a_0',a_1',\cdots, a_n')&=&I(a_0\alpha^n, a_1\alpha^{n-1}\delta, \cdots,a_n\delta^n)\nonumber\\
&=&\sum
\{Z_{\nu_0\nu_1\cdots\nu_n}{a_0}^{\nu_0}\alpha^{n\nu_0}{a_1}
^{\nu_1}\}\alpha^{(n-1)\nu_1}\delta^{\nu_1}\nonumber\\ &&\cdots
{a_i}^{\nu_i}\alpha^{(n-i)\nu_i}\delta^
{i\nu_i}\cdots{a_n}^{\nu_n}\delta^{n\nu_n}\nonumber\\
&=&\sum \{Z{\nu_0\nu_1\cdots\nu_n}{a_0}^{\nu_0}{a_1}^{\nu_1}\cdots
{a_i}^{\nu _i}\cdots {a_n}^{\nu_n}\}\nonumber\\
&&
.\alpha^{n\nu_0+(n-1)\nu_1+\cdots+(n-i)\nu_i+\cdots+\nu_{n-1}}\delta^{\nu_1+
\cdots+i\nu_i+n\nu_n}\nonumber\\
&=&\alpha^p\beta^p\sum
Z_{\nu_0\nu_1\cdots\nu_n}{a_0}^{\nu_0}{a_1}^{\nu_1}\cdots
{a_n}^{\nu_n},\nonumber
\end{eqnarray}
from which we obtain the identities
\begin{eqnarray}
n\nu_0+(n-1)\nu_1+\cdots+(n-i)\nu_i+\cdots+\nu_{n-1}=p,\nonumber\\
\nu_1+\cdots+i\nu_i+\cdots+(n-1)\nu_{n-1}+n\nu_n=p.\nonumber
\end{eqnarray}
Addition of this formula results in:
$n(\nu_0+\nu_1+\cdots+\nu_n)=2p$. \hfill\ $\Box$

\medskip So, two characteristic equations we have found are the
following
\begin{eqnarray}
\nu_1+2\nu_2+3\nu_3+\cdots+n\nu_n=p,\qquad ng=2p.\nonumber
\end{eqnarray}
On the other hand, the theorem just proven has evidently a
converse.
\paragraph{Theorem 2.5}
{\it Every homogeneous function of the coefficients a, which satisfies the equation
$ng-2p=0,$ where $g$ is the degree, $p$ the weight of this
function, and $n$ the order of the base form, is an invariant with
respect to transformation $(5)$.}

\medskip The proof is evident; if one substitutes the $a'$ instead
of the $a$ into $I$, then a factor $\alpha^p \delta^p$ appears
because of the assumed identities. So, we avoid repeating it.
\section {The Operation Symbols $\mathcal{D}$ and $\Delta$}
We introduce two special operational symbols
\begin{eqnarray}
\mathcal{D}&=&a_0\frac{\partial}{\partial
a_1}+2a_1\frac{\partial}{\partial
a_2}+3a_2\frac{\partial}{\partial a_3}+\cdots
+na_{n-1}\frac{\partial}{\partial a_n},\nonumber\\
\Delta &=& n a_1\frac{\partial}{\partial a_0}+(n-1)a_2
\frac{\partial}{\partial a_1}+(n-2)a_3\frac{\partial}{\partial
a_2}+\cdots+a_n\frac{\partial}{\partial a_{n-1}}.\nonumber
\end{eqnarray}
that are two infinitesimal generators of the unimodular group.
\paragraph{Lemma 3.1}
 {\it The operations $\mathcal{D}$ and
$\Delta$ applied to a homogeneous function, leave the degree $g$
in the $a_i$ unchanged; $\mathcal{D}$ lowers the weight $p$ by 1,
that is, to $p-1$; $\Delta$ raises it by 1, to $p+1$. Also,
homogeneity is preserved by these operations.}
\paragraph{Theorem 3.2}
{\it If $\mathcal{A}$ is a homogeneous function in the $a_i$, of degree $g$ and weight $p$, then
\begin{eqnarray}
({\mathcal{D}} \Delta - \Delta {\mathcal{D}}) {\mathcal{A}} =(n
g-2p){\mathcal{A}}
\end{eqnarray}}

\medskip\noindent{\it Proof: } We proof this theorem in three steps:

1. If formulae $(9)$ is valid for two two expressions
${\mathcal{A}}_1$ and ${\mathcal{A}}_2$ which have the same degree
$g$ and the same weight $p$, then it is also valid for the sum
${\mathcal{A}}_1+{\mathcal{A}}_2.$ From
\begin{eqnarray}
({\mathcal{D}} \Delta -\Delta {\mathcal{D}}){\mathcal{A}}_1&=&(n
g-2p){\mathcal{A}}_1,\nonumber\\
({\mathcal{D}} \Delta -\Delta {\mathcal{D}}){\mathcal{A}}_2&=&(n
g-2p){\mathcal{A}}_2\nonumber
\end{eqnarray}
it follows immediately by addition that
\begin{eqnarray}
({\mathcal{D}} \Delta -\Delta
{\mathcal{D}})({\mathcal{A}}_1+{\mathcal{A}}_2)=(n
g-2p)({\mathcal{A}}_1+{\mathcal{A}}_2)\nonumber
\end{eqnarray}

2. If the formulae is valid for ${\mathcal{A}}_1$(degree$=g_1,$
weight$=p_1$) and ${\mathcal{A}}_2$(degree$=g_2,$ weight$=p_2$),
then it is also valid for the product ${\mathcal{A}}_1
{\mathcal{A}}_2$. The product is easily seen to be a homogeneous
function of degree $g=g_1+g_2$ and weight $p=p_1+p_2$. According
of our assumptions, we have now
\begin{eqnarray}
({\mathcal{D}} \Delta - \Delta {\mathcal{D}}){\mathcal{A}}_1=(n
g_1-2p_1){\mathcal{A}}_1,\nonumber\\
({\mathcal{D}} \Delta - \Delta {\mathcal{D}}){\mathcal{A}}_2=(n
g_2-2p_2){\mathcal{A}}_2.\nonumber
\end{eqnarray}
Using the rules, one finds for the product:
\begin{eqnarray}
({\mathcal{D}} \Delta - \Delta {\mathcal{D}})({\mathcal{A}}_1
{\mathcal{A}}_2)
&=&{\mathcal{D}} \Delta ({\mathcal{A}}_1{\mathcal{A}}_2)-\Delta
{\mathcal{D}}({\mathcal{A}}_1{\mathcal{A}}_2)\nonumber\\
&=&{\mathcal{D}}\{{\mathcal{A}}_1\Delta
{\mathcal{A}}_2+{\mathcal{A}}_2\Delta{\mathcal{A}}_1\}-\Delta
\{{\mathcal{A}}_1{\mathcal{D}}{\mathcal{A}}_2+{\mathcal{A}}_2{\mathcal{D}}{\mathcal{A}}_1\}\nonumber\\
&=&{\mathcal{A}}_1{\mathcal{D}}\Delta {\mathcal{A}}_2+\{\Delta
{\mathcal{A}}_2\}{\mathcal{D}}{\mathcal{A}}_1\}+{\mathcal{A}}_2{\mathcal{D}}\Delta
{\mathcal{A}}_1+\{\Delta
{\mathcal{A}}_1\}\{{\mathcal{D}}{\mathcal{A}}_2\}\nonumber\\
&&-{\mathcal{A}}_1\Delta
{\mathcal{D}}{\mathcal{A}}_2-\{{\mathcal{D}}{\mathcal{A}}_2\}\{\Delta
{\mathcal{A}}_1\}-{\mathcal{A}}_2\Delta
{\mathcal{D}}{\mathcal{A}}_1-\{{\mathcal{D}}{\mathcal{A}}_1\}\{\Delta
{\mathcal{A}}_2\}\nonumber\\
&=&{\mathcal{A}}_1\{{\mathcal{D}}\Delta {\mathcal{A}}_2-\Delta
{\mathcal{D}}{\mathcal{A}}_2\}+{\mathcal{A}}_2\{{\mathcal{D}}\Delta{\mathcal{A}}_1-\Delta
{\mathcal{D}}{\mathcal{A}}_1\}\nonumber\\
&=&(n g_2-2p_2){\mathcal{A}}_1{\mathcal{A}}_2+(n g_1
-2p_1){\mathcal{A}}_2{\mathcal{A}}_1\nonumber\\
&=&(n(g_1+g_2)-2(p_1+p_2)){\mathcal{A}}_1{\mathcal{A}}_2;\nonumber
\end{eqnarray}
Therefore, indeed
\begin{eqnarray}
({\mathcal{D}} \Delta - \Delta
{\mathcal{D}})({\mathcal{A}}_1{\mathcal{A}}_2)=(n
g-2p){\mathcal{A}}_1 {\mathcal{A}}_2\nonumber
\end{eqnarray}

3. The formulae is valid for ${\mathcal{A}}=Const$, since then it
only says $0=0$, and for ${\mathcal{A}}=a_i$, we have
\begin{eqnarray}
{\mathcal{D}}\Delta a_i&=&{\mathcal{D}}(n-i)a_{i+1}=(n-i)(i+1)a_i,\nonumber\\
\Delta{\mathcal{D}}a_i&=&\Delta i a_{i-1}=i(n-i+1)a_i.\nonumber
\end{eqnarray}
Consequently,
\begin{eqnarray}
({\mathcal{D}}\Delta-\Delta{\mathcal{D}}a_i)=(n
i-i^2+n-i-in+i^2-i)a_i =(n-2i)a_i.\nonumber
\end{eqnarray}
Finally the general theorem follows directly from these three
facts.\hfill\
$\Box$

\medskip We shall need to derive two general formulas from formulae $(9)$.
Let ${\mathcal{A}}$ be a homogeneous function of degree $g$ and
weight $p$. Then ${\mathcal{D}}{\mathcal{A}}$ is a homogeneous
function of degree $g$ and weight $p-1$, and generally, if we let
$${\mathcal{D}}^2={\mathcal{D}}{\mathcal{D}},\qquad
{\mathcal{D}}^3={\mathcal{D}}{\mathcal{D}}{\mathcal{D}},\qquad
\cdots$$ then ${\mathcal{D}}^k{\mathcal{A}}$ is a homogeneous
function of degree $g$ and weight $p-k.$ Thus, with
straightforward computations we have the following formulas :
\begin{eqnarray}
{\mathcal{D}}^2 \Delta - \Delta {\mathcal{D}}^2&=&2(n
g-2p+1){\mathcal{D}},\nonumber\\
{\mathcal{D}}\Delta^2 - \Delta^2{\mathcal{D}}&=&2(n
g-2p-1)\Delta.\nonumber
\end{eqnarray}
\paragraph{Theorem 3.3}
 {\it The following two important formulas are
valid in general:
\begin{eqnarray}
{\mathcal{D}} \Delta - \Delta {\mathcal{D}}=k(n g
-2p+k-1){\mathcal{D}}^{k-1}, \\
{\mathcal{D}} \Delta ^k - \Delta ^k {\mathcal{D}}=k(n g
-2p-k+1)\Delta^{k-1}.
\end{eqnarray}}

\medskip\noindent{\it Proof : } This is easily proven by induction from $k$
to $k+1$. Namely, suppose the formulas $(10)$, $(11)$  are valid
for $k$. We first apply the operation ${\mathcal{D}}$ to $(10)$,
then formula $(9)$ to ${\mathcal{D}}^k$, likewise, we first apply
the operation $\Delta$ to $(11)$, then formulae $(9)$ to $\Delta
^k$. Then we obtain the following identities :
\begin{eqnarray}
{\mathcal{D}}^k+1 \Delta - {\mathcal{D}}\Delta
{\mathcal{D}}^k&=&k(n g -2p+k-1){\mathcal{D}}^k,\nonumber\\
{\mathcal{D}}\Delta {\mathcal{D}}^k-\Delta {\mathcal{D}}^k+1&=&(n
g-2(p-k)){\mathcal{D}}^k,\nonumber\\
\Delta {\mathcal{D}}\Delta ^k- \Delta^{k+1}{\mathcal{D}}&=&k(n
g-2p-k+1)\Delta^k,\nonumber\\
{\mathcal{D}}\Delta^{k+1}-\Delta{\mathcal{D}}\Delta^k&=&(n
g-2)(p+k))\Delta^k.\nonumber
\end{eqnarray}
And if one adds the first two formulas on the one hand and the
second two on the other hand, then one obtains formulas $(10)$ and
$(11)$ for $k+1$. But since they are valid for $k=2$, they are
valid in general. For $k=1$, they both turn into formulae $(9)$,
provided one defines
${\mathcal{D}}^0{\mathcal{A}}=\Delta^0{\mathcal{A}}={\mathcal{A}}$
\hfill\ $\Box$

We shall assume that the formulas we derived are always
being applied to a homogeneous function.
\section{The smallest system of conditions for the determination of the invariants}
Using the relations derived in the previous section, we can
substantially simplify the necessary and sufficient conditions
for invariants.
\paragraph{Theorem 4.1}
{\it Every homogeneous function $I$ of the
coefficients $a_0,a_1,\cdots ,a_n$ of degree $g$ and weight $p$,
where $n g=2p$, is an invariant if it satisfies the differential
equation ${\mathcal{D}}I=0$.}

\medskip\noindent{\it Proof: } Let ${\mathcal{A}}$ be a homogeneous
function of degree $g$ and weight $p$, where $n g=2p$. Suppose
further that ${\mathcal{A}}$ satisfies the differential equation
${\mathcal{D}}{\mathcal{A}}=0$. Then we want to show that
${\mathcal{A}}$ also satisfies the differential equation $\Delta
{\mathcal{A}}=0$, and hence is an invariant.

Since $\Delta ^k {\mathcal{A}}$ has degree $g$ for all $k$, the
weight can not take on arbitrary values; therefore, from some
$\Delta ^i$ on, all have to be identically zero. Let this be
$\Delta^k{\mathcal{A}}=0$, while $\Delta^{k-1}{\mathcal{A}}\neq
0$. But, according to formulae $(III)$  in the preceding section,
we have
\begin{eqnarray}
{\mathcal{D}}\Delta
{\mathcal{A}}-\Delta^k{\mathcal{D}}{\mathcal{A}}=k(-k+1)\Delta^{k-1}{\mathcal{A}}\nonumber
\end{eqnarray}
Therefore, we must have $k(k-1)\Delta^{k-1}{\mathcal{A}}=0$, Here
it can neither happen that $\Delta^{k-1}{\mathcal{A}}=0$ (because
of our assumption), nor that $k=0$, because then ${\mathcal{A}}$
would be identically zero, and so it follows that $k=1$, that is
$\Delta{\mathcal{A}}=0$, and our assertion is proven. \hfill\
 $\Box$

\end{document}